\documentclass{article}

\usepackage[english]{babel}

\usepackage[letterpaper,top=2cm,bottom=2cm,left=3cm,right=3cm,marginparwidth=1.75cm]{geometry}

\usepackage{amsmath}
\usepackage{amssymb}
\usepackage{mathrsfs}
\usepackage{graphicx}
\usepackage{xspace}
\usepackage[colorlinks=true, allcolors=blue]{hyperref}

\title{The classification of the finite groups whose abelian subgroups of equal prime power order are conjugate}
\author{Robert W. van der Waall}

\begin{document}
\maketitle
\begin{flushleft}
\begin{abstract}
Let $G$ be a finite group and assume $p$ is a prime dividing the order of $G$. Suppose for any such $p$, that every two abelian $p$-subgroups of $G$ of equal order are conjugate. The structure of such a group $G$ has been settled in this article.
\end{abstract}

\textbf{MSC-2020 Classifications}\\
20D06, 20D08, 20D15, 20D20, 20D25, 20D30, 20D40, 20E07, 20E45, 20F16, 20F28.\\

\textbf{Key Words and Phrases}\\
Finite groups, conjugacy classes, generalized Fitting subgroup, Frattini subgroup, abelian $p$-subgroup, nilpotent groups, supersolvable groups, simple groups of Lie type, alternating and symmetric groups, sporadic simple groups, automorphism groups, Schur multiplier.

\tableofcontents

\section{Introduction}

Let us consider the following classes of finite groups.
\newline
\newline
$(1) \text{  }\mathscr{B}=\{G|\text{ each pair of subgroups of $G$ of equal order consists of conjugate subgroups} \}$
\newline
\newline
$(2)  \text{  }\mathscr{H}=\{G| \text{ each pair of supersolvable subgroups of $G$ of equal order consists of conjugate}$ \\
$\text{subgroups} \}.$
\newline
\newline
In the formulation of $(2)$, the word \textit{supersolvable} can successively be replaced by \emph{nilpotent}, \emph{abelian}, \emph{cyclic}. In doing so, one gets the corresponding classes of finite groups $\mathscr{N}, \mathscr{A},\mathscr{C}$. Next, let us change the phrase \emph{of equal order} from $(1)$ and $(2)$ into \emph{of equal prime power order}; as such, one gets the classes $(1)': \mathscr{B}_{\pi}$ and $(2)': \mathscr{H}_{\pi}$ respectively. In addition, one gets the classes of finite groups $\mathscr{N}_{\pi}, \mathscr{A}_{\pi},\mathscr{C}_{\pi}$ by replacing in $(2)'$ the word \emph{supersolvable} again successively by \emph{nilpotent}, \emph{abelian}, \emph{cyclic}.
\newline
\newline
The classification of the groups in the classes $\mathscr{B}, \mathscr{N}, \mathscr{H}$ and $\mathscr{A}$ has been subject of study in the papers $[3],[8],[9],[10]$ and $[11]$; that of $\mathscr{C}$ can be read off from $[2]$. Notice that the following hierarchies of classes of groups are apparent: $\mathscr{B} = \mathscr{H} \subsetneq \mathscr{N} \subsetneq \mathscr{A} \subsetneq \mathscr{C}$ (the former $=$-sign is in fact Theorem $2.2$ in $[11]$), and $\mathscr{B}_{\pi} \subseteq \mathscr{H}_{\pi} \subseteq \mathscr{N}_{\pi} \subseteq \mathscr{A}_{\pi} \subseteq \mathscr{C}_{\pi}$. Since any group of prime power order is nilpotent, any abelian group is nilpotent and any nilpotent group is supersolvable, this forces immediately $\mathscr{B}_{\pi} = \mathscr{H}_{\pi} = \mathscr{N}_{\pi}$. As it will be obtained below, it holds that the class of groups $\mathscr{A}_{\pi}$ is properly contained in $\mathscr{C}_{\pi}$, and likewise $\mathscr{N}_{\pi} \subsetneq \mathscr{A}_{\pi}$. The simple first Janko group $J_1$ is a member of the class $\mathscr{B}_{\pi}$, but not of $\mathscr{B}$; see $[11]$. An example of a solvable group in the class $\mathscr{B}_{\pi}$, but not appearing in the class $\mathscr{B}$, has been found by F. Gross; see pages $337$ and $338$ in $[4]$. Hence $\mathscr{B} \subsetneq \mathscr{B}_{\pi}$ anyway.
\newline
\newline
In this paper all groups are finite. Notations, symbols etc., will be standard; see for instance $[5]$, but also $[3]$.
\newline
In $[7]$ the structure of the solvable groups from the class $\mathscr{C}_{\pi}$ has essentially been determined, while in $[12]$ the structure of the non-solvable groups belonging to $\mathscr{C}_{\pi}$ has been established, next to a plethora of other results.
\newline
\newline
We will exhibit the structure of the groups belonging to the classes $\mathscr{B}_{\pi}$ and $\mathscr{A}_{\pi}$. Salient results from $[7]$ and $[12]$ will be properly quoted when needed.

\section{Non-solvable $\mathscr{A}_{\pi}$-groups}
We start with the classification of the non-solvable $\mathscr{A}_{\pi}$-groups, i.e. non-solvable $G \in \mathscr{A}_{\pi}$. Let $X \in \mathscr{C}_{\pi}$, then we know from $(2.3)$ Theorem, $[12]$ that any chief factor of $X$ is either cyclic of prime order, or non-cyclic elementary abelian, or, it is isomorphic to a specific non-abelian simple group $N$. In the last situation, $N$ happens to be isomorphic to one of the simple groups from the following list $\mathscr{L}$, consisting of the groups labeled $a)$ up to and including $h)$, see $\mathsection 5$ in $[12]$.
\newline
\newline 
$a)$ the Mathieu simple groups $M_{11}$ and $M_{23}$; 
\newline
$b)$ the simple first Janko group $J_1$;
\newline
$c)$ the simple groups $PSL(2,q)$, $q=p^n$, $p$ an odd prime, $n \geq 1$, $n$ odd, $q \geq 7$;
\newline
$d)$ the groups $SL(2,q)$, $q=p^n$, $p$ an odd prime, $n \geq 1$, $n$ odd, $q \geq 7$;
\newline
$e)$ the simple Suzuki group $Sz(2^{2m+1}), m \geq 1$;
\newline
$f)$ the simple group $PSL(2,2^a), a \geq 3$;
\newline
$g)$ $SL(2,5)$;
\newline
$h)$ the simple alternating group $A_5$.
\newline
\newline
It has been shown in $\mathsection 5$ of $[12]$ that any group from the list $\mathscr{L}$ belongs to $\mathscr{C}_{\pi}$. It will be implicitly used that no two simple groups in the list $\mathscr{L}$ are isomorphic to each other, see $[1]$. In addition, it has been shown in $(6.1)$ Theorem, $[12]$, that in any $X \in \mathscr{C}_{\pi}$ admitting a non-solvable chief factor $L/K (K,L \unlhd X, K \leq L)$, it is true that any odd prime dividing $|L/K|$ does not divide the product $|X/L| \cdot |K|$; either one of the groups $X/L$ and $K$ happens to be solvable, see $(2.3)$ Theorem, $[12]$.
 
\subsection{Case a)}
In applying $\mathscr{A}_{\pi} \subseteq \mathscr{C}_{\pi}$, one observes that none of the groups $M_{11}$ or $M_{23}$ can be isomorphic to a chief factor of some alleged $G \in \mathscr{A}_{\pi}$, in its role of $X \in \mathscr{C}_{\pi}$ mentioned above. [Indeed, suppose it were, then $(4.2)$ Theorem, $[12]$ guarantees the existence of an $R \unlhd G$ with $G \cong M \times R$, $M \in \{M_{11},M_{23}\}$ and $(|M|,|R|)=1$ . Thus $G \in \mathscr{A}_{\pi}$ would imply $M \in \mathscr{A}_{\pi}$. Now, as each Sylow $2$-subgroup of $M_{11}$ and of $M_{23}$ is neither cyclic, nor elementary abelian, it follows on the contrary from the $\mathscr{A}_{\pi}$-property for $G$ when applied to the abelian $2$-subgroups of $G$ of equal order being conjugate, that the alleged statement "$M \in \mathscr{A}_{\pi}$" in untenable. All this forces that the earlier appropriate group $N$ is not a chief factor of $G \in \mathscr{A}_{\pi}$ isomorphic to each of $M_{11}$ and $M_{23}$. Hence such a group $G$ is not an $\mathscr{A}_{\pi}$-group.]

\subsection{Case e)}
Now assume one deals with a group $X \in \mathscr{C}_{\pi}$ in which a chief factor might be isomorphic to some  $Sz(2^{2m+1}), m \geq 1$. Then $(5.5)$ and $(5.6)$ Theorems, $[12]$ reveal that a normal subgroup $S$ of $X$ exists with $S \cong Sz(2^{2m+1})$. Any Sylow $2$-subgroup of $S$ does contain cyclic subgroups of order $4$, as well as subgroups isomorphic to $C_2 \times C_2$. Therefore, the group $X \in \mathscr{C}_{\pi}$ is not an $\mathscr{A}_{\pi}$-group.
\newline
\newline
Next, we are going to deal successively with the classes of groups $b), c), d), f)$ followed by $g)$ and $h)$.

\subsection{Case b)}
Due to $(4.2)$ Theorem, $[12]$, any $X \in \mathscr{C}_{\pi}$ admitting a chief factor isomorphic to the first Janko group $J_1$, does satisfy the structure $X \cong N \times R$, in which $N \cong J_1$ and where $(|N|,|R|)=1$, with $R \unlhd X, N \unlhd X, [N,R]=1$. It holds that $J_1 \in \mathscr{N}$, see Theorem $2.1$, $[11]$. As surely now $J_1 \in  \mathscr{A}_{\pi}$ holds, the specialized assumption $X \in \mathscr{A}_{\pi}$ implies though that also $R \in \mathscr{A}_{\pi}$; notice that $R$ is solvable as $2 \nmid |R|$. The contents of Theorem C in $[7]$ imply that the Sylow subgroups of $R$ are cyclic or non-cyclic elementary abelian. Hence as $(|J_1|,|R|)=1$ and $X \in \mathscr{A}_{\pi}$ (whence $X \in \mathscr{C}_{\pi}$) we conclude that in particular $X \in \mathscr{B}_{\pi}$. Indeed, $J_1 \in \mathscr{N}_{\pi}=\mathscr{B}_{\pi}$ and $R \in \mathscr{A}_{\pi}$ implies here $R \in \mathscr{B}_{\pi}$.

\subsection{Case c)}
Let $X \in \mathscr{C}_{\pi}$ admit an isomorphic copy of a non-abelian simple group occurring in $c)$, as a chief factor $L/K$ and assume the order of $K$ is odd. Then it is known from $(5.9)$ Theorem, $[12]$, that there exists an $N \unlhd X$, with $N \cong L/K$ and $(|X/N|,|N|)=1$. Any Sylow subgroup of $X/N$ is cyclic or elementary abelian, see Theorem C in $[7]$. The Sylow $2$-subgroups of the groups from the list $\mathscr{L}$ occurring in $c)$ are dihedral of order $2^h, h \geq 3$ (say) or elementary abelian of order $4$. The Sylow $t$-subgroups, $t \neq 2$ of $N$ are cyclic or non-cyclic elementary abelian. 
\newline
Next, we specialize these conditions to $X \in \mathscr{A}_{\pi}$. Hence, by the order formula $(p^n-1) \cdot p^n \cdot (p^n+1)=|PSL(2,p^n)|=|N|$, where $N \cong PSL(2,p^n), N \unlhd X, X \in \mathscr{A}_{\pi}$, $N \in \mathscr{L}$, one deduces that the Sylow $2$-subgroups of $X$ and hence of $N$ must be isomorphic to $C_2 \times C_2$. It follows that $X$ actually belongs to $\mathscr{B}_{\pi}$! In addition, the aforementioned order formula provides the number theoretic condition $p \equiv \pm 3$ mod $8$.

\subsection{Case d)}
It has been shown in $(2.2)$ Theorem in $[7]$ that any factor group of a group $X \in \mathscr{C}_{\pi}$ belongs to $\mathscr{C}_{\pi}$. Assume that $K \unlhd X$ exists in which $\tilde{L} \unlhd X$ with $\tilde{L} > K$ satisfies that $\tilde{L}/K$ is isomorphic to a simple group from $c)$ in the list $\mathscr{L}$. Assume $K$ is of even order. Then $(5.9)$ Theorem, $[12]$ yields the existence of a normal subgroup $S$ of $X$ with $S \cong SL(2,p^n)$, with this $SL(2,p^n)$ being a group in $d)$ from the list $\mathscr{L}$. Notice that here $(|X/S|.|S|)=1$ by that same Theorem $(5.9)$, that $X/S \in \mathscr{C}_{\pi}$ and that also $\tilde{L}/K$ is isomorphic to $SL(2,p^n)/Z(SL(2,p^n))$, i.e. isomorphic to $S/Z(S)$. Remember $X/\tilde{L}, K$ and $X/S$ are all solvable by the Jordan-Hölder-Zassenhaus Theorem. By the way, the logically equivalent statement "$SL(2,t) \in \mathscr{C}_{\pi} \iff PSL(2,t) \in \mathscr{C}_{\pi}$", where $t$ is some specific power of an odd prime, is certainly non-trivial to figure out: one has to use the facts that the Sylow $2$-subgroups of $SL(2,t)$ are generalized quaternion, whereas the corresponding Sylow $2$-subgroups of $PSL(2,t)$ are dihedral or isomorphic to $C_2 \times C_2$. Now, if one specializes to $X \in \mathscr{A}_{\pi}$ in this Case $d)$, then certainly $X \cong SL(2,q) \rtimes R$ for some specific $q$ being an odd power of an odd prime, $(|R|,|SL(2,q)|)=1$, the Sylow subgroups of $R$ being cyclic or non-cyclic elementary abelian. All the non-cyclic elementary abelian Sylow subgroups of $R$ do centralize $SL(2,q)$, following $(5.8)$ Theorem, $[12]$. Hence by Sylow, the implication $X \in \mathscr{A}_{\pi} \Rightarrow R \in \mathscr{B}_{\pi}$ holds true. It is possible under the assumption $X \in \mathscr{A}_{\pi}$ to pin down more specific information on the group $SL(2,q)$. At first, $SL(2,q) \in \mathscr{C}_{\pi}$ by $(5.9)$ Theorem, $[12]$, implying $PSL(2,q) \in \mathscr{C}_{\pi}$. Furthermore, all abelian $2$-subgroups of $SL(2,q)$ are in fact cyclic, since a Sylow $2$-subgroup of $SL(2,q)$ is generalized quaternion. It is true, but a little bit tricky to verify, that any two cyclic subgroups of equal order of $SL(2,q)$ are conjugate in $SL(2,q)$, see for instance $[11]$. Now let $q=p^n$, $p$ an odd prime, $n \geq 1$ odd. Notice $X/C_X(S) \hookrightarrow Aut(SL(2,p^n))$ and that $Out(SL(2,p^n))$ is isomorphic to a cyclic group of order $2n$, just as $n$ is odd. Moreover, due to $(5.4)$ Theorem, $[12]$, it holds for $X \in \mathscr{A}_{\pi} \subseteq \mathscr{C}_{\pi}$, that $X$ is isomorphic to $(SL(2,p^n) \times U)\langle \varphi \rangle$, with $|\varphi|$ dividing $n$, that $(|U\langle \varphi \rangle|,|SL(2,p^n)|)=1$ and that $U \leq C_{S \rtimes U \langle \varphi \rangle}(S)$. The Sylow $l$-subgroups of $X$ with $l$ any prime dividing $|U\langle \varphi \rangle/U|$ are cyclic, whereas $U\langle \varphi \rangle$ has precisely cyclic or non-cyclic elementary abelian Sylow subgroups, due to $(5.4)$ Theorem, $[12]$. Hence, since $X \in \mathscr{A}_{\pi}$. it follows that $U\langle \varphi \rangle \in \mathscr{A}_{\pi}$.
\newline 
We did show that the abelian $2$-subgroups of our $X \in \mathscr{A}_{\pi}$ are contained in $S \unlhd X$ and that the Sylow $u$-subgroups coincide with the Sylow $u$-subgroups of $X$ for any prime $u$ dividing $|S|$. We are going to show that the integer $n$ occurring in $SL(2,p^n) \cong S$ turns out to be equal to $1$ or to $3$. This runs as follows.
\newline
Assume $n \geq 3$. The Sylow $p$-subgroups of $S$, whence of $X \in \mathscr{A}_{\pi}$, are non-cyclic elementary abelian. The Frattini-argument (I$.7.8$ Satz, $[5]$) provides $X=SN_X(S_p)$, where $S_p \in Syl_p(S)=Syl_p(X)$. Therefore, $X/S$, being isomorphic to $N_X(S_p)/N_S(S_p)$, is solvable, as $ 2 \mid |S|$ and $(|X/S|,|S|)=1$. Notice that also $N_{S/Z(S)}((S_pZ(S)/Z(S))$ is isomorphic to $N_S(S_p)Z(S)/Z(S)$, whence to $N_S(S_p)/Z(S)$. The group $N_{S/Z(S)}((S_pZ(S)/Z(S))$ is solvable, as $S_pZ(S)/Z(S)$ is elementary abelian of order $p^n$ and as $N_{S/Z(S)}((S_pZ(S)/Z(S))/(S_pZ(S)/Z(S))$ is cyclic; see II$.8.21$ Hilfsatz $[5]$. Therefore, each of the groups $N_X(S_p)/N_S(S_p)$, $N_S(S_p)/Z(S)$ and $Z(S)$ is solvable, whence the group $N_X(S_p)$ itself is solvable. Now, remember, $S_p$ is elementary abelian, $S_p \in Syl_p(X)$ and $X \in \mathscr{A}_{\pi}$. Thus $X$ permutes all of the $p$-subgroups of equal order under conjugation action. In combination with $N_X(S_p)$ being solvable, one observes that Lemma $3$, $[4]$, due to Gross, is usable, yielding $n = 3$. 
\newline
The proof of Lemma $3$, $[4]$ is a quite a little bit involved, among others, it needs Huppert's Transitivity Theorem for doubly-transitive solvable groups; see XII.$7.3$ Theorem, $[5b]$. Originally in this rubric Case d), the integer $n$ in the group $SL(2,p^n)$ satisfies $n \equiv 1$ (mod $2$). Hence indeed the integer $n$ equals $1$ or $3$. Thus $S$ is isomorphic to $SL(2,p)$ or to $SL(2,p^3)$. Observe, such an $S$ belongs $\emph{a fortiori}$ to $\mathscr{A}$; see Main Theorem, $[8]$. Note also that now $X$ is isomorphic to a direct product of the group $SL(2,p^{\delta})$ (with $\delta=1$ or $\delta=3$) and a group $R$ satisfying  $(|SL(2,p^{\delta})|,|R|)=1, R \in \mathscr{B}_{\pi}$.
\newline 
[A little warning. Only when for an odd prime $p$, $ p \equiv \pm 3$ (mod $8$) happens to be the case, it is true that $SL(2,p)$ and $SL(2,p^3)$ are groups living in  $\mathscr{B}_{\pi}$. Indeed, otherwise $16$ would  divide the order of a Sylow $2$-subgroup $\overline{Q}$  of $SL(2,p)$ and likewise of $SL(2,p^3)$, yielding that $\overline{Q}$ would contain one cyclic subgroup of order $|\overline{Q}|/2$ and simultaneously two subgroups of that order $|\overline{Q}|/2$ being (generalized) quaternion.]

\subsection{Case f)}
We infer from $(5.5)$ Theorem, $[12]$, that there exists an $N \unlhd Y$ with $(|N|,|Y/N|)=1$ and $N \cong PSL(2,2^a)$, whenever it is assumed that $Y \in \mathscr{C}_{\pi}$ and that $Y$ admits a chief factor $L/K$ being isomorphic to the group $PSL(2,2^a)$ in case $a \geq 3$. The Sylow $2$-subgroups of $N$ are non-cyclic elementary abelian. The other Sylow subgroups of $N$ are cyclic. Now, as the odd primes dividing $|N|$ are distinct from those dividing $|Y/N|$, it appears that any Sylow $t$-subgroup of $Y/N$ happens to be cyclic or non-cyclic elementary abelian; indeed, $t \geq 5$ holds, so that $(5.8)$ Theorem, $[12]$ yields the conclusion. If we specialize to $Y \in \mathscr{A}_{\pi}$ here, one deduces though that $Y$ is henceforth included in $\mathscr{B}_{\pi}$. 
\newline 
\newline
In summary, we did show in the cases $b), c)$ and $f)$, that the assumption $G \in \mathscr{A}_{\pi}$ may lead to $G \in \mathscr{B}_{\pi}$. Precisely said, in each of these three cases any non-abelian simple normal subgroup $N$ of $G$ originating from the cases $b), c)$ and $f)$ with $G/N$ solvable and satisfying $(|N|,|G/N|) = 1$ with $G \in \mathscr{A}_{\pi}$, does yield $G \in \mathscr{B}_{\pi}$, up to the situation in which
$p \equiv \pm 1$ (mod $8$) holds with respect to the group $N \cong PSL(2,p^n)$. 
\newline
\newline
There is more to say about those cases $b), c)$ and $f)$. In $[4]$ it has been brought to light that any factor group $T/E$ with $E \unlhd T$ and $T \in \mathscr{B}_{\pi}$, belongs to $\mathscr{B}_{\pi}$. Furthermore, due to 
Theorem 1, $[4]$, it holds that, given a non-solvable group $T \in \mathscr{B}_{\pi}$, the factor group $T/O_{2'}(T)$ is isomorphic to one of the groups from the list $\mathscr{V}$, defined by $\mathscr{V} = \{J_1, PSL(2,p), PSL(2,p^3), SL(2,8), SL(2,p), SL(2,p^3), P \mathit{\Gamma} L(2,32)\}$, where $p$ is any prime congruent to $\pm 3$ (mod $8$) and where  $P \mathit{\Gamma}L(2,32)$ stands for the group $PSL(2,32) \rtimes \langle \vartheta \rangle$; here $\vartheta \in Aut(PSL(2,32))$ with $|\vartheta| = 5$. The automorphism $\vartheta$ operates like a Frobenius automorphism on the elements inside any matrix of the group $SL(2,32)$; notice $SL(2,32) = PSL(2,32)$. As Gross did observe, each group from $\mathscr{V}$ belongs to $\mathscr{B}_{\pi}$. Of course, the solvable group $O_{2'}(T)$ stands for the normal subgroup of $T$ subject to being maximal of odd order.  
\newline
\newline
The results obtained so far give rise to the following theorems.
\newline
\newline
$\textbf{Theorem 1}$ Let $G \in \mathscr{A}_{\pi}$. Assume there exists a non-abelian chief factor of $G$, isomorphic to $PSL(2,2^a)$ for some $a \geq 3$. Then either $a=5$, or else $G \cong PSL(2,8) \times R$, where $(|R|,|PSL(2,8)|)=1$, and $R \in \mathscr{B}_{\pi}$. It holds that $PSL(2,8) \in \mathscr{B}_{\pi}$; in fact, $PSL(2,8) \in \mathscr{B}$.
\newline
\newline
\textbf{Proof} The assumptions of the theorem have been considered above. Those did lead to the implication $G \in \mathscr{A}_{\pi} \Rightarrow G \in \mathscr{B}_{\pi}$. Thus in consulting the list $\mathscr{V}$ and applying the Jordan-Hölder-Zassenhaus Theorem, we see that either $a=3$ or $a=5$ is forced. One observes that the direct product property holds for $a=3$; use Case $f)$ and Theorem $1$, $[4]$, with factors $N \cong PSL(2,8)$ and  $R \unlhd G$, generating $G$ and satisfying $(|N|,|R|) = 1$. Groups in $\mathscr{B}_{\pi}$ are factor-group-closed in taking homomorphisms. Thus $R \in \mathscr{B}_{\pi}$ as $R \cong G/N$.   The fact that $PSL(2,8)$ belongs to $\mathscr{B}$ can be found in $[3]$. $\square$   
\newline
\newline
Similar to the proof of Theorem $1$, one is able to show the truth of the statements in Theorem $2$.
\newline
\newline
$\textbf{Theorem 2}$  Let $G \in \mathscr{A}_{\pi}$. Assume $L/K$ is a chief  factor of $G$ being non-abelian. Suppose $L/K$ is neither isomorphic to $A_5$ nor isomorphic to $PSL(2,32)$. In addition, suppose $|K|$ is odd. Then $G$ is a direct product of a group $N$ and a group $R$, where $N$ stems from the list $\mathscr{V}$ as to the cases $b)$ and $c)$ of the list $\mathscr{L}$, such that $(|N|,|R|) = 1$ holds. Notice $G \in \mathscr{B}_{\pi}$ and $R \in \mathscr{B}_{\pi}$. $\square$  
\newline
\newline
More involved is the proof of the next Theorem 3, in which one deals with a group $G \in \mathscr{A}_{\pi}$ admitting a chief section isomorphic to $PSL(2,32)$.  
\newline
\newline
$\textbf{Theorem 3}$ Let $G \in \mathscr{A}_{\pi}$ and assume $G$ admits a chief factor containing a subgroup isomorphic to $PSL(2,32)$. Then $G =(N \times B)\langle \alpha \rangle$, satisfying $N \unlhd G$, $N \cong PSL(2,32)$, $B \unlhd G$, $(|N|,|B\langle \alpha \rangle|) = 1$, $\alpha \in G \setminus(N \times B)$ some $5$-element of $G$ with $\alpha^5 \in B$ and with $\langle \alpha \rangle \in Syl_5(G)$. In addition, $G \in \mathscr{B}_{\pi}$ is fulfilled and  $B\langle \alpha \rangle \in \mathscr{B}_{\pi}$ too. Furthermore $N\langle \alpha \rangle \in \mathscr{B}$. 
\newline
\newline
$\textbf{Proof}$ By $(5.3), (5.5), (6.1)$ Theorems,$[12]$, one gets that there exists $N \unlhd G$ satisfying $(|N|,|G/N|)= 1$ and $N \cong PSL(2,32)$. The group $N$ is simple. One has $G/C_G(N) \hookrightarrow Aut(N)$, where $Aut(N) \cong N \rtimes \langle \vartheta \rangle$ with $|\vartheta|=5$. Look at $NC_G(N)/C_G(N)$. That group is isomorphic to $N/(N \cap C_G(N))$, whence to $N$. Since $G \in \mathscr{A}_{\pi}$, it cannot be that $G$ equals $NC_G(N)$, which is $N \times C_G(N)$. [Indeed, suppose the contrary. Thus one has $(|N|,|C_G(N)|) = 1$ where also $N \in \mathscr{A}_{\pi}$ would happen. However, it would follow that $N \in \mathscr{B}$, because the Sylow $2$-subgroups of $N$ are elementary abelian and the other Sylow subgroups of $N$ cyclic. It is known, see page $400$, $[3]$, that $PSL(2,32)$ is $\emph{not}$ a member of $\mathscr{B}$. So we did produce a contradiction].  Therefore, due to Theorem $1$, $[4]$ in combination with the Jordan-Hölder-Zassenhaus Theorem, one gets 
$G/C_G(N) \cong (PSL(2,32) \rtimes \langle \vartheta \rangle)$. Thus $C_G(N)$ is the wanted group $B \unlhd G$, showing there exists a $5$-element $\alpha \in G\setminus (N \times B)$ satisfying $G = (N \times B)\langle \alpha \rangle$. The Sylow $t$-subgroups of $B$ are elementary abelian or cyclic for $t$ any prime with $t \geq 7$ and $t$ dividing $|B|$. Subgroups of $B$ of equal order of some power of the prime $t$ are conjugate by means of products of elements of $B$ and powers of $\alpha$, whereas the Sylow $5$-subgroups of $B\langle \alpha \rangle$ are cyclic. [Indeed, any Sylow $5$-subgroup of $G$ is conjugate to a Sylow $5$-subgroup of $B \langle \alpha \rangle$; it holds that $G \in \mathscr{A}_{\pi}$ (given), so that $\alpha^5 \in B$. In case $\alpha^5 \neq 1$ this forces that Sylow $5$-subgroups of $B$ cannot be non-cyclic elementary abelian; if not, any Sylow $5$-subgroup of $G$ would be non-cyclic elementary abelian and contained in $B$ itself, which is not the case.] The conclusion is that $B\langle \alpha \rangle \in \mathscr{B}_{\pi}$. As to  $N\langle \alpha \rangle \in \mathscr{B}$, see Theorem $9.1$, $[3]$. Notice $G \in \mathscr{B}_{\pi}$ follows. $\square$
\newline 
\newline  
$\textbf{Theorem 4}$ Let $G \in \mathscr{A}_{\pi}$ and assume $N \unlhd G$ exists with $N \cong SL(2,5)$ or $N \cong A_5$. Then $G = N \times R$ with $R \unlhd G$ satisfying $(|N|,|R|) = 1$. Furthermore, $SL(2,5) \in \mathscr{B}, A_5 \in \mathscr{B}$ and $R \in \mathscr{B}_{\pi}$, implying $G \in \mathscr{B}_{\pi}$. 
\newline
\newline
$\textbf{Proof}$ Do combine the outcome of $(5.4)$ Theorem, $[12]$ with that of Theorem $1$, $[4]$. The result is the direct product property, yielding $R \in \mathscr{B}_{\pi}$. The statements $SL(2,5) \in \mathscr{B}$ and $A_5 \in \mathscr{B}$ are to be found in Lemma $6$ and Theorem $11$, $[3]$. $\square$ 
\newline
\newline
In accordance with the results proved so far, one observes that details of the structure of a non-solvable group $G \in \mathscr{A}_{\pi}$, are related to those groups $U \in \mathscr{A}_{\pi}$ satisfying $(6,|U|)=1$. Those groups $U$ are solvable by the Feit-Thompson Theorem. So let us assume $N$ is a non-cyclic minimal normal subgroup of $U$. Hence $N$ is elementary abelian of $p$-power order for some prime $p \geq 5$. Since $U \in \mathscr{A}_{\pi}$, $Syl_p(U)$ is a singleton, i.e. $\{N\}= Syl_p (U)$, due to $(4.4)$ Theorem, $[7]$. Look at the group $O_{p'}(U)$, being the normal subgroup of $U$ of highest possible order relatively prime to $p$. Now, any non-cyclic Sylow subgroup of our $U \in \mathscr{A}_{\pi}$ with $(6,|U|) = 1$, is elementary abelian and normal in $U$; see the proof of Theorem $12$ later on. Therefore $O _{p'}(U)N/O_{p'}(U)$ is not only contained in the Fitting subgroup $F(U/O_{p'}(U))$ of $U/O_{p'}(U)$, but even better, as $Syl_p(O_{p'}(U)N/O_{p'}(U)) = \{O_{p'}(U)N/O_{p'}(U)\} \cong \{N\}$ and as $F(U/O _{p'}(U))$ is nilpotent, satisfying 
$Syl_p(U/O_{p'}(U))=Syl_p(O_{p'}(U)N/O_{p'}(U))$, it holds that $F(U/O_{p'}(U))=O_{p'}(U)N/O_{p'}(U)$.  Above it was implicitly observed that all non-cyclic Sylow $t$-subgroups of $U$ (if any) are contained in $O_{p'}(U)$, whenever such prime $t$ exists unequal to $p$. Hence all Sylow subgroups of $U/O_{p'}(U)$ are cyclic, whence $U/O_{p'}(U)$ is meta-cyclic; see Theorem $5$ and $(6.2)$ Theorem, $[7]$. Hence the group $U/(O_{p'}(U)N)$, being isomorphic to $(U/O_{p'}(U))/F(U/O_{p'}(U))$, does permute all subgroups of equal order of the elementary abelian group $O _{p'}(U)N/O_{p'}(U)$ transitively under conjugation action, as $U \in \mathscr{A}_{\pi}$. So we have been ended up into the contents of $(7.2)$ Theorem, $[7]$, whence it follows that $U/O_{p'}(U) \in \mathscr{B}$ holds. Therefore, due to Theorem $7.2$, $[10]$), the order of $O_{p'}(U)N/O_{p'}(U)$ equals $p^3$; here $(6,|U|)=1$ is used again. Also, via Theorem $5.3$ of $[10]$, one gets that   
$U/(O_{p'}(U)N) \cong (U/O_{p'}(U))/F(U/O_{p'}(U))$  implies, using Theorem $7.2$, $[10]$, that $U/(O_{p'}(U)N)$ is cyclic. Thus we have got that for any prime $p_i$ dividing $|U|$ for which  $S_i \in Syl_{p_i}(U)$ is not cyclic, it is true that $S_i$ is a normal subgroup of $U$ and that $S_i$ is elementary abelian of order $p_i^3$; whence, if there are $d \geq 2$ such primes $p_i$, one has    $\prod \limits_{j=1}^d |S_j| = \prod \limits_{j=1}^d p_j^3$  and also $U/\hat{S_i} \in \mathscr{B}$, where $\hat{S_i}$ stands for the product $S_1 \cdot S_2 \cdots \check{S_i} \cdots S_d$, meaning that $S_i$ is $\emph{not}$ included in taking the product of the $S$'s. Finally, when it happens that it occurs that $(p_r^2 + p_r + 1, p_s^2 + p_s + 1) = 1$ whenever $ 1 \leq r < s \leq d$, then $U \in \mathscr{B}$ holds; see Corollary $5.5$ of Theorem $5.3$, $[10]$. Conversely, if $U \in \mathscr{A}_{\pi}$ is given with $(6,|U|) = 1$, then $U \in \mathscr{B}$, provided that relatively-prime-condition regarding the $p_j$'s $(j = 1, \cdots, d)$ from the last sentence is due to hold; see again Theorem $5.3$, $[10]$. Furthermore, if the previous integer $d$ happens to be equal to $1$, then $U$ appears to be a member of $\mathscr{B}$; see Theorem $3$, $[4]$. Notice that any $U \in \mathscr{A}_{\pi}$ with $(6,|U|) = 1$ is a member of the class $\mathscr{B}_{\pi}$.  
\subsection{Cases g) and h)}
In order to handle the cases g) and h) from the list $\mathscr{L}$, we will prove the following Theorem $\emph{inter alia}$, irrespective of the structure of any arbitrary $G \in \mathscr{A}_{\pi}$.
\newline 
\newline
$\textbf{Theorem 5}$ Let $G \in \mathscr{A}_{\pi}$. Assume  $N \unlhd G$ exists with  $2 \nmid |N|$. Then $G/N \in \mathscr{A}_{\pi}$.            
\newline
\newline
$\textbf{Proof}$ Suppose $G$ is a minimal counterexample to the Theorem. Then there exists a normal subgroup $T$ of $G$ for which $G/T \notin \mathscr{A}_{\pi}$ and $|T|$ odd. Take $T$ $\emph{a fortiori}$ as small as possible with respect to that conclusion. Then $T$ is a minimal normal subgroup of $G$. [Indeed, otherwise there exists a $V \unlhd G$ with $T \gneqq V \gneqq     \{1\}$, here any such $V$ satisfies  $G/V \in \mathscr{A}_{\pi}$. The group $V$ is of odd order. Also $G/T \cong (G/V)/(T/V)$, whence, as $V \neq  \{1\}$, $G/T \in \mathscr{A}_{\pi}$ does follow, contrary to the choice of $T$ in $G$.] Now, due to $|T|$ being odd, the group $T$ is solvable, whence $T$ is a $p$-group being cyclic or elementary abelian for some specific odd prime $p$. Let  $X/T$ and $Y/T$ be abelian subgroups of $G/T$  of equal prime power order,  $G \geq X$ and $G \geq Y$.  Assume firstly that  $p$  does not divide $|X/T|$ . Hence  by the Schur-Zassenhaus Theorem, there exists  $X_1 \leq G$ and $Y_1 \leq G$ of equal order satisfying  $X=TX_1, Y=TY_1$  and  $X_1 \cap T= \{1\} = Y_1 \cap T$. Moreover, as we know, the abelian groups $X_1$ and $Y_1$ are of equal prime power order. Since  $G \in \mathscr{A}_{\pi}$, the groups $X_1$ and $Y_1$ are conjugate in $G$. Hence $X/T$ and $Y/T$ are conjugate in $G/T$ and of equal prime power order. Thus we may also assume that $X/T$ and $Y/T$ are $p$-groups. Suppose that $T$ were non-cyclic elementary abelian. Then, as $G \in \mathscr{A}_{\pi}$, no Sylow $p$-subgroup $P$ of $G$ does contain a cyclic subgroup of order $p^2$. Hence, any Sylow $p$-subgroup of $X$ (and of $Y$) would be contained in $T$, due to the $\mathscr{A}_{\pi}$-property of $G$. This, however, contradicts $p$ dividing $|X/T|$ with $|X/T| = |Y/T|$. Therefore, because $T$ is a minimal normal subgroup of $G$, it follows that $T$ is cyclic of order $p$. Thus, as $p$ is odd, Huppert's Theorem III.$8.3$ Satz, $[5]$, forces the Sylow $p$-subgroups of $G$ to be cyclic, since $G \in \mathscr{A}_{\pi}$ yields $T$ being the only existing subgroup of order $p$ in G.  As $X/T$ and $Y/T$ are $p$-subgroups of $G/T$ of equal order, one observes that $X$ and $Y$ are $p$-subgroups of $G$ of equal order (remember $X/T$ and $Y/T$ were chosen to be of equal prime power order). Now  $X$ is contained in some cyclic Sylow $p$-subgroup $S$ of $G$ and $Y$ is likewise contained in some cyclic Sylow  $p$-subgroup $\overline{S}$ of $G$. It holds by Sylow Theory, that $S$ and $\overline{S}$ are conjugate in $G$. Therefore in particular $X$ and $Y$ are conjugate in $G$. Then, however, the cyclic groups $X/T$ and $Y/T$, being of equal prime power order, are conjugate in $G/T$. All in all now, one gets that, irrespective of the particular prime power order of $X/T$ (and of $Y/T$), it always holds that $X/T$ and $Y/T$ are conjugate in $G/T$. This, however, contradicts the choice of $T$ appearing in the second sentence of this proof with respect to $G/T \notin \mathscr{A}_{\pi}$. Therefore, no such $T$ exists. The Theorem has been proved. $\square$    
\newline
\newline
$\textbf{Remark}$ The condition $2 \nmid |N|$ in Theorem $5$ cannot be omitted. The group $SL(2,7)$ belongs to $\mathscr{A}$, whence to $\mathscr{A}_{\pi}$, but $PSL(2,7)$ being a factor group of $SL(2,7)$, does contain cyclic subgroups of order $4$ as well as non-cyclic abelian subgroups of order $4$. Therefore $PSL(2,7)$ is $\emph{not}$ contained in $\mathscr{A}_{\pi}$. 
\newline 
\newline 
$\textbf{Theorem 6}$ Let $G \in \mathscr{A}_{\pi}$    and assume there exists a chief factor $L/K$ of $G$ containing a subgroup isomorphic to $A_5$. Then $L/K \cong A_5$. Suppose $2 \nmid |K|$. Then there exists an $N \unlhd G$  with  $N \cong A_5$. It holds that $G =N  \times  R$, a direct product of the group $N$ and some group $R \in \mathscr{B}_{\pi}$ with $(|N|,|R|) = 1$. 
\newline
\newline
$\textbf{Proof}$ Look at $(5.8)$ Theorem, $[12]$ and the fact that $R \in \mathscr{A}_{\pi}$ with $(6,|R|) = 1$ implies that $R \in \mathscr{B}_{\pi}$, as we saw above. $\square$     
\newline
\newline
In the next theorem we quote $(7.1)$ Theorem, $[12]$, in connection to Theorem $4$ in shorthand shape. 
\newline
\newline
$\textbf{Theorem 7}$ Suppose $L/K \cong A_5$ is a chief factor of $G \in \mathscr{C}_{\pi}$. Assume $2 \mid |K|$. Then
\newline
$1)$ $\emph{either}$ $G \cong SL(2,5) \times U$,  with $U \in \mathscr{C}_{\pi}$ and $(30,|U|) = 1$;
\newline
$2)$ $\emph{or else}$ $G = (V \times T) \rtimes W, V \unlhd G, T \unlhd G, (|V|,|T|)=1=(|V \times T|,|W|)$; 
$W$ is a meta-cyclic $\mathscr{C}_{\pi}$-group, or $W=\{1\}$, or $W$ is cyclic; $T=\{1\}$ or $T=\prod $(non-cyclic elementary abelian Sylow subgroups $S$ of $TW$) with all $S \unlhd TW$; $TW$ is a solvable $\mathscr{C}_{\pi}$-group, $[[W,W],F]=\{1\}, F \lhd G, V=F \rtimes S, [V,V]=V, S \cong SL(2,5)$; $F=\prod \limits_{i=1}^4 (C_{p_i} \times  C_{p_i})^{\delta_i}$, with $p_1=11, p_2=19, p_3=29, p_4=59$, all $\delta_i \in \{0,1\}$ except $\delta_1=\delta_2=\delta_3=\delta_4=0.$ $\square$
\newline
\newline
Next, in Theorem $8$, we do specialize the group $G$  occurring in Theorem $7$, in being a group from $\mathscr{A}_{\pi}$.
\newline
\newline
$\textbf{Theorem 8}$ Let $G$ be a group from the assertion in Theorem $7$ and assume  $G \in \mathscr{A}_{\pi}$. Then, using the notations and notions from Theorem $7$,
\newline 
$1)$ $\emph{either}$ $G \cong SL(2,5) \times U$,  with $U \in \mathscr{B}_{\pi}$ and $(30,|U|) = 1$;
\newline
$2)$ $\emph{or else}$  $G/F \in \mathscr{A}_{\pi},  G/(F \rtimes Z(S)) \in \mathscr{B}_{\pi}, [S,TW] = \{1\}$;
\newline
$TW$ is a solvable $\mathscr{B}_{\pi}$-group all of whose non-cyclic Sylow $t$-subgroups are of order $t^3$ and also elementary abelian, whenever $t$ is a prime dividing $|T|$;
\newline
$W \in \mathscr{B}$, and if $W \neq \{1\}$, then $W$ is (meta-)cyclic satisfying $(|W/[W,W]|,|[W,W]|)=1,  VW \in \mathscr{A}_{\pi}$ and $[[W,W],F] = \{1\}$. 
\newline
\newline
$\textbf{Proof}$ $\emph{Re 1)}$ 
Here $G \in \mathscr{A}_{\pi}$ has certainly the structure of part $1)$ in the announcement of Theorem $7$. In addition, $U \in \mathscr{B}_{\pi}$ holds due to the results discussed earlier between Theorem $4$ and Theorem $5$.
\newline
$\emph{Re 2)}$ Because of Theorem $6$, $G/F \in \mathscr{A}_{\pi}$ holds. Hence $G/F$ satisfies a structure as found in part $1)$ of Theorem $7$ and also that from Re $1)$ just proven regarding Theorem $8$. Therefore also $G/(F \rtimes Z(S))$ is an $\mathscr{A}_{\pi}$- group; notice $PSL(2,5) \in \mathscr{B}$ and $V \in \mathscr{B}_{\pi}$, together with 
$G/(F \rtimes Z(S)) \cong A_5 \times TW$. The rest of the Theorem has been discussed earlier. The property $W\in \mathscr{B}$ has been shown in Theorem $3$, $[4]$ and its immediate Corollary; notice that it is implicitly used that each of the Sylow subgroups of $W$ is equal to some Sylow subgroup of $G$ and that those subgroups are cyclic. Furthermore $VW \cong G/T \in \mathscr{A}_{\pi}$ holds; the groups $[W,W]$ and $F$ do centralize each other, as shown in 
$(7.1)$ Theorem, $[12]$. $\square$
\newline
\newline
$\textbf{Remark}$ It is surely possible that in the Theorems $7$ and $8$ the group $W$ does not act trivially on $F$ by conjugation. It might only happen when  $p_3=29$ with $\delta_3=1$, where $W$ contains a cyclic Sylow $7$-subgroup $\langle c \rangle$ acting non-trivially on $F$ properly by conjugation action such that $[ \langle c^7 \rangle ,F]= \{1\}$ and when $p_4=59$ with $\delta_4= 1$ occurs and $W$ possesses a cyclic Sylow $29$-subgroup $\langle d \rangle$ such that $[\langle d^{29} \rangle, F]=\{1\}$. Notice, however, inside the corresponding group $G \in \mathscr{A}_{\pi}$, it is not possible that both the described situations happen to occur simultaneously. The easiest paradigm examples of Theorem $8$, part $2)$, are 
$G\cong (C_{29} \times  C_{29}) \rtimes     (SL(2,5)  \times  C_7)$  and $G \cong (C_{59} \rtimes  C_{59}) \rtimes (SL(2,5) \times C_{29})$; both these groups belong to $\mathscr{B}$, whence to $\mathscr{B}_{\pi}$. See Theorem $11$, $[3]$ and $\mathsection$$7$, $[12]$.
\subsection{Sylow subgroups of non-solvable  $\mathscr{A}_{\pi}$-groups}
We are going to find out, what structure a Sylow subgroup of a given $Y \in \mathscr{A}_{\pi}$ can  afford. Notice that a non-cyclic $t$-subgroup of an non-solvable $Y \in \mathscr{A}_{\pi}$, $t$ an odd prime, is indeed elementary abelian, see  the next Theorem $9$. As to that statement regarding solvable $\mathscr{A}_{\pi}$-groups, see Section $3$ below. 
\newline
\newline
$\textbf{Theorem 9}$ Let $Y \in \mathscr{A}_{\pi}$ be non-solvable and assume $S_p \in Syl_p(Y)$, $p$ an odd prime. Then $S_p$ is either cyclic or else it is non-cyclic elementary abelian. 
\newline
\newline
$\textbf{Proof}$ Suppose firstly $p \geq 5$ and assume $S_p$ is not cyclic. Then, by $(5.8)$ Corollary, $[12]$, $S_p$ is a normal subgroup of $Y$ or some $PSL(2,p^f) (f \geq 2)$ is isomorphic to a chief section of $Y$. If $S_p \unlhd Y \in \mathscr{A}_{\pi}$ happens to be true, then $S_p$  turns out to be elementary abelian. Thus suppose that some $PSL(2,p^f) (f \geq 2)$ is isomorphic to a chief factor of $Y$, whence $f$ is odd. Hence $S_p$ is elementary abelian because of $S_p$  being isomorphic to a Sylow $p$-subgroup of that $PSL(2,p^f)$; see  II.$8.10$ Satz, $[5]$ and $(6.1)$ Theorem, $[12]$.  What about $p=3$? Assume $S_3$  is not cyclic and that $S_3$ and the Fitting subgroup $F(Y)$ of $Y$ do intersect each other non-trivially. Then it follows that $S_3$ is in fact a normal subgroup of $Y$, due to $Y \in \mathscr{A}_{\pi}$ in combination with III.$7.5$ Hilfssatz, $[5]$, so that $S_3$ is elementary abelian too, by the fact that here $Y \in \mathscr{A}_{\pi}$ is $3$-solvable implying, by a result of Shult VIII.$7.11$ Remarks, $[5b]$, that $S_3$ happens to be abelian at first sight; see also $(4.1)$ Lemma and Theorem C, $[7]$. Thus assume $S_3 \cap F(Y)=\{1\}$. It follows then from the classification of the non-solvable $\mathscr{C}_{\pi}$-groups as carried out in $(5.4)$ , $(6.1)$ and $(7.1)$ Theorems, $[12]$, that either $S_3$ is isomorphic to a subgroup of a simple chief factor of $ \in \mathscr{A}_{\pi}$, in which indeed our non-cyclic group $S_3$ turns out to be elementary abelian, or else that $Sz(2^{2m+1}) \leq Y \in \mathscr{A}_{\pi}$ should be fulfilled. The last eventuality however, is not possible, as the groups $Sz(2^{2m+1})$ do possess cyclic as well as non-cyclic subgroups all being of order $4$. Therefore, for any $S_3 \in Syl_3(Y)$ with $Y \in \mathscr{A}_{\pi}$, the assumption on $S_3$ being non-cyclic does lead to the fact that $S_3$ turns out to be elementary abelian. $\square$
\newline
\newline
Let us introduce some notation. $E_{p^n}$ --  an elementary abelian group of order $p^n$ ($p$ any prime, $n$ a natural number); $Q_{\nu}$ -- a generalized quaternion group of order $2^{\nu}$ with  $\nu \geq 4$; $Q$ -- the quaternion group of order $8$; $C_{u}$ -- a cyclic group of order $u$ with $u \geq 1$.
\newline
\newline
$\textbf{Theorem 10}$ Let $X \in \mathscr{A}_{\pi}$ be non-solvable. Then $S_2 \in Syl_2(X)$ is isomorphic to one of the following groups:  $Q_{\nu}, Q, E_4, E_8,  E_{32}$. 
\newline
\newline
$\textbf{Proof}$ If $S_2$ would be cyclic, then $X$ would be solvable, as it is known from Burnside's theorem that $X/O_{2'}(X) \cong S_2$ (see IV.$2.8$ Satz, $[5]$); remember, the odd order group $O_{2'}(X)$ is solvable due to the Feit-Thompson Theorem. If $S_2$ is not cyclic, then the statement regarding the structure of $S_2$ is to be found in the results sofar. $\square$

\section{Solvable $\mathscr{A}_{\pi}$-groups}
In this Section we do focus our attention on the structure of the Sylow subgroups of solvable $\mathscr{A}_{\pi}$-groups. To start with, let us look at the following Theorem.
\newline
\newline
$\textbf{Theorem 11}$ Let $X \in \mathscr{C}_{\pi}$ be solvable and assume $p$ is a prime dividing $|X|$. Then $S_p \in Syl_p(X)$ is isomorphic to a group belonging to the following classes of groups: $\{C_{p^a}| a \geq 1\}$, $\{E_{p^b}|b \geq 2\}$, $\{Q\}$, $\{\text{Suzuki } 2\text{-groups} \neq Q\}$.
\newline
\newline
$\textbf{Proof}$ The statement of the Theorem is the outcome of $\mathsection$$3$, $[7]$. $\square$
\newline
\newline
Let us specialize Theorem $11$ to the case in which $X \in \mathscr{A}_{\pi}$ is solvable. Notice, that if a group $Y \in \mathscr{C}_{\pi}$ admits Suzuki $2$-groups $\neq Q$ as Sylow $2$-subgroups, it cannot be that $Y$ belongs to $\mathscr{A}_{\pi}$; as such, those Suzuki $2$-groups do contain cyclic subgroups of order $4$ as well as non-cyclic abelian subgroups of order $4$.
\newline
\newline
$\textbf{Theorem 12}$ Let $X \in \mathscr{A}_{\pi}$ be solvable and assume $p$ is a prime dividing $|X|$. Then $S_p \in Syl_p(X)$ is isomorphic to a group from the following classes of groups: $\{C_{p^a}|a \geq 1\}$, $\{E_{p^2}|p \neq 2\}$,$\{E_{p^3}|p \neq 2\}$, $\{E_4\}$, $\{E_8\}$, $\{E_{32}\}$, $\{Q\}$. Furthermore, $X/O_{2'}(X)$ is isomorphic to one of the following groups:  $C_{2^a} (a \geq 0), E_4 \rtimes C_3, E_8 \rtimes (C_7 \rtimes C_3), E_8 \rtimes C_7, E_{32} \rtimes (C_{31} \rtimes C_5), Q \rtimes C_3$.
\newline
\newline
$\textbf{Proof}$ By using the property that for any given prime $t$, Sylow $t$-subgroups of any group are conjugate to each other in that group, it follows from Theorem $11$, that given $X \in \mathscr{A}_{\pi}$ being solvable and observing that $Q$ has only as subgroups of order $4$ the cyclic ones and that it contains only one subgroup of order $2$, it holds that $X \in \mathscr{B}_{\pi}$. [This result will be mentioned as Theorem $15$, in order to give it a prominent place.]. Hence by Theorem $1$, $[4]$) one observes that any Sylow $2$-subgroup of $X$ is either isomorphic to some cyclic $2$-group, or to $E_4$, or to $E_8$, or to $E_{32}$, or to $Q$. Namely, Gross in his Theorem $1$, $[4]$ did show the truth of the statement as formulated in the sentence of our Theorem beginning with 
"Furthermore  ... ". We proceed in providing self-contained elaborated reasoning around the structure of $S_p \in Syl_p(X)$, when $p$ is odd. Suppose $S_p$ is not cyclic. Then there do not exist cyclic subgroups of order $p^i$   in $S_p$ for any $i \geq 2$.  [Indeed, assume it were. Then, as $S_p \leq X \in \mathscr{A}_{\pi}$, $S_p$ contains exactly one subgroup of order $p$. [Otherwise, as $Z(S) \geq C_p$ for some specific $C_p$, there would exist $C_p \times D$ with $D \neq C_p, D \leq S_p, |D| = p$; a contradiction to $X \in \mathscr{A}_{\pi}$.] Hence III.$7.5$ Hilfssatz, $[5]$ yields $S_p$ being cyclic, a contradiction to the assumed non-cyclicity of $S_p$.] Thus $S_p$ does contain only as non-trivial cyclic subgroups, its cyclic subgroups of order $p$. Now  $X \in \mathscr{A}_{\pi}$ permutes all its cyclic subgroups of prime power order transitively under conjugation action. Therefore, by 
VIII.$5.8.b$ Theorem and VIII.$7.11$ Remarks, $[5a]$ in which a lemma of Shult is involved, one concludes that $S_p$ is abelian. Thus, as $S_p$ does not contain elements of order $p^2$, $S_p$ is elementary abelian. Hence Gross' Lemma $3$, $[4]$ is valid for our non-cyclic group $S_p$  with $p \neq 2$ too, i.e.  $S_p \cong E_{p^2}$ or $S_p \cong E_{p^3}$. $\square$ 
\newline
\newline
$\textbf{Corollary 13}$ Assume $S_p \in Syl_p(X)$ happens to be non-cyclic for some prime $p$, where $X \in \mathscr{A}_{\pi}$ is understood to be solvable. Then $S_p \unlhd X$ or $S_p \cong Q$.
\newline
\newline
$\textbf{Proof}$ Apply Theorem $12$. Hence one gets, by using $(4.4)$ Theorem, $[7]$, that  $S_p \unlhd X$ unless $S_p \cong Q$. [Indeed, $(4.4)$ Theorem, $[7]$ says that elementary abelian non-cyclic Sylow subgroups of a solvable group $Y \in \mathscr{C}_{\pi}$, happen to be entirely contained in the Fitting subgroup of $Y$, whence implying these subgroups to be normal in $Y$.] $\square$
\newline
\newline
$\textbf{Corollary 14}$ Assume $Q$ to be isomorphic to some Sylow $2$-subgroup $S_2$ of some $X \in \mathscr{A}_{\pi}$. Then not always $S_2 \unlhd X$ does hold.
\newline
\newline
$\textbf{Proof}$ Take $X \cong E_{5^2} \rtimes SL(2,3)$, where the action by conjugation of $SL(2,3)$ on $E_{5^2}$ is Frobenius and faithful. Then $X \in \mathscr{B} \subseteq \mathscr{A}_{\pi}$, but $S_2 \in Syl_2(X)$ being isomorphic to $Q$, is not a normal subgroup of $X$. Notice $X$ is a solvable group. $\square$
\newline
\newline
To be complete, when $S_2$ is a cyclic Sylow $2$-subgroup of some group $H$, it holds that $H/O_{2'}(H) \cong S_2$; see IV.$2.8$ Satz, $[5]$.
\newline
\newline
$\textbf{Theorem 15}$ Let $G \in \mathscr{A}_{\pi}$ be solvable, then $G \in \mathscr{B}_{\pi}$.
\newline
\newline
$\textbf{Proof}$ The truth of the statement of this theorem follows directly from combining Theorem $12$ and Corollary $13$. $\square$   
\newline
\newline
Due to Theorem $15$, one observes that solvable $\mathscr{A}_{\pi}$-groups are closed under taking homomorphisms, as follows.
\newline
\newline
$\textbf{Theorem 16}$ Let $N$ be a normal subgroup of a solvable group $G \in \mathscr{A}_{\pi}$. Then $G/N \in \mathscr{A}_{\pi}$.
\newline
\newline
$\textbf{Proof}$ It follows from Theorem $15$ that  $G \in \mathscr{B}_{\pi}$ holds. Groups in $\mathscr{B}_{\pi}$ are closed under taking homomorphisms, as observed on page $333$, $[4]$. Thus $G/N \in \mathscr{B}_{\pi}$. As $\mathscr{B}_{\pi} \subseteq \mathscr{A}_{\pi}$, we are done. $\square$
\newline
\newline
There exists a theorem, a converse to Theorem $16$ and valid for any $G$, solvable or not, as follows.
\newline
\newline
$\textbf{Theorem 17}$ Let $G$ be a group and assume $N_1 \unlhd G$ and  $N_2 \unlhd G$ do exist, satisfying  
$(|N_1|,|N_2|)=1$. Suppose $G/N_1 \in \mathscr{A}_{\pi}$ and $G/N_2 \in \mathscr{A}_{\pi}$. Then $G \in \mathscr{A}_{\pi}$.
\newline
\newline
$\textbf{Proof}$ Let $X$ and $Y$ be abelian subgroups of $G$ of equal prime-power order  $p^m (m \geq 1)$. Then at least one of the orders $|N_1|$ or $|N_2|$, say $|N_1|$, is not divisible by $p$. Hence $XN_1/N_1$, being isomorphic to $X/(X \cap N_1 )$, and $YN_1/N_1$, being isomorphic to
$Y/(Y \cap N_1)$, are abelian groups each of order $p^m$, just by  $X \cap N_1=\{1\}=Y \cap N_1$. Therefore, as $G/N_1 \in \mathscr{A}_{\pi}$, there exists $g \in G$,  satisfying $(gN_1)^{-1}(XN_1/N_1)(gN_1)=YN_1/N_1$. Thus $(g^{-1}Xg)N_1/N_1=YN_1/N_1$. The groups $g^{-1}Xg$ and $Y$ are Sylow $p$-subgroups of $YN_1$. Hence, by Sylow's theorem, there exists a $t \in YN_1$, such that  $Y =t^{-1}(g^{-1}Xg)t=(gt)^{-1}X(gt)$. In other words, $X$ and $Y$ are conjugate in $G$. So indeed, $G \in \mathscr{A}_{\pi}$ does follow. $\square$ 
\newline
\newline
Summarizing, according to Theorem 15, any solvable  $\mathscr{A}_{\pi}$-group has turned out to be 
a $\mathscr{B}_{\pi}$-group. while it holds that any non-cyclic Sylow $p$-subgroup $S \not \cong Q$ of a solvable $\mathscr{A}_{\pi}$-group $G$ is normal in $G$; namely those groups $S$ are elementary abelian, making $(4.4)$ Theorem, $[7]$ in vogue. Due to these facts, the analysis of the structure of a solvable $\mathscr{A}_{\pi}$-group is rather easy to accomplish. Let us see what happens.
\newline
\newline
$\textbf{Theorem 18}$ Let $G \in \mathscr{A}_{\pi}$ be solvable and assume there exists $E_{p^2} < G$  for some particular prime $p$. Then
\newline
$1)$ $\emph{either}$ $G/O_{p'}(G)$ does contain a normal Sylow $p$-subgroup being non-cyclic elementary abelian, whereas the other Sylow subgroups of $G/O_{p'}(G)$  (if any) are cyclic;
\newline
$2)$ $\emph{or else}$ $p=5$ and  $G/O_{5'}(G) \cong E_{5^2} \rtimes SL(2,3)$, or $p=11$ and $G/O_{11'}(G) \cong E_{11^2} \rtimes SL(2,3)$, or $G/O_{11'}(G)  \cong E_{11^2} \rtimes (SL(2,3) \times C_5)$. 
\newline
\newline
$\textbf{Proof}$ Notice $G/O_{p'}(G)$ is solvable and $G/O_{p'}(G) \in \mathscr{A}_{\pi}$. We want to apply the contents of $(7.2)$ Theorem, $[7]$ with respect to the group $G/O_{p'}(G)$. In order to accomplish that goal, one has to do the following.
\newline
$\emph{Assume}$ $G/O_{p'}(G)$ does contain a non-cyclic non-elementary Sylow subgroup; whence it has to be non-abelian due to Theorem $12$, for the  group $G/O_{p'}(G) \in \mathscr{A}_{\pi}$. It holds that  $E_{p^2} \leq S_p \in Syl_p(G)$, so that $S_p$ is elementary abelian. By $(4.4)$ Theorem, $[7]$ it follows that $S_p \lhd G$. Thus $S_pO_{p'}(G)/O_{p'}(G)$, being a Sylow $p$-subgroup of $G/O_{p'}(G)$, happens to be normal in $G/O_{p'}(G)$. So $p$ does not divide $|G/(S_pO_{p'}(G))|$. By definition of $O_{p'}(G)$, the fact that $S_pO_{p'}(G)/O_{p'}(G)$ is normal in the Fitting subgroup $F(G/O_{p'}(G))$, the solvability of $G/O_{p'}(G)$ and the fact that $F(G/O_{p'}(G))$ is nilpotent, all together imply  that $F(G/O_{p'}(G))=S_pO_{p'}(G)/O_{p'}(G)$. Notice now, as $G/O_{p'}(G) \in \mathscr{A}_{\pi}$, that $S_pO_{p'}(G)/O_{p'}(G)$ is a chief factor of $G/O_{p'}(G)$ and that there exists no other minimal normal subgroup in $G/O_{p'}(G)$ than  $S_pO_{p'}(G)/O_{p'}(G)$, due to the solvability of
$G/O_{p'}(G)$ and $F(G/O_{p'}(G))$ being the unique normal Sylow $p$-subgroup of $G/O_{p'}(G)$. Hence 
$G/O_{p'}(G)$ is a so-called $(P)$-subdirectly irreducible solvable $\mathscr{A}_{\pi}$-group in the terminology of $[7]$, containing a non-abelian Sylow subgroup by our assumption in the beginning of the proof of the Theorem. Hence we have reached the point where $(7.2)$ Theorem, $[7]$ can be applied  
on $G/O_{p'}(G)$, namely that it enforces that $p$ divides the integer $55$; as such one gets the structure of the three particular groups appearing in the statement of the Theorem. The working-out of those facts has been done in the proof of $(7.2)$ Theorem, $[7]$.
\newline
$\emph{Next suppose}$ that all Sylow subgroups of $G/O_{p'}(G) \in \mathscr{A}_{\pi}$ are abelian. One of them is the non-cyclic elementary abelian and normal $p$-subgroup  $S_pO_{p'}(G)/O_{p'}(G)$.  Any other alleged proper Sylow subgroup $\overline{S}$  of $G/O_{p'}(G) \in \mathscr{A}_{\pi}$ is cyclic or otherwise non-cyclic elementary abelian because of  $G/O_{p'}(G) \in \mathscr{A}_{\pi}$. If $\overline{S}$ is not cyclic, it thus would be normal in $G/O_{p'}(G)$ by Corollary $13$. However, this is not possible, for then 
$\overline{S} \leq F(G/O_{p'}(G))=$ $ S_pO_{p'}(G)/O_{p'}(G)$ produces the contradiction
$p \nmid |\overline{S}| \mid |S_pO_{p'}(G)/O_{p'}(G)|=|S_p |$. Hence all Sylow subgroups of $G/O_{p'}(G)$ other than $S_pO_{p'}(G)/O_{p'}(G)$ are cyclic. The Theorem has been proved. $\square$
\newline
\newline
$\textbf{Corollary 19}$ Let $G \in \mathscr{A}_{\pi}$ be solvable. Assume there exists a prime $p$ dividing $|G|$ such that $E_{p^2} < G$. Then $G/O_{p'}(G) \in \mathscr{B}$.
\newline
\newline
$\textbf{Proof}$ Look at the statements in Theorem $18$. If $G/O_{p'}(G)$ is a group occurring in the phrase just after the word "either", then $G/O_{p'}(G) \in \mathscr{B}$ holds as it satisfies the contents of Theorem $11$, $[9]$. As to the three exceptional groups remaining in Theorem $18$, a detailed study by means of the computer language GAP has been carried out in Sektion $3.2$, $[6]$. $\square$
\newline
\newline
All the facts about solvable $\mathscr{A}_{\pi}$-groups obtained sofar, can be collected in the following Portmanteau Theorem.
\newline
\newline
$\textbf{Theorem 20}$ (Portmanteau Theorem)
\newline
Let $G$ be a solvable $\mathscr{A}_{\pi}$-group. Then one of the following eight statements is true and all do occur in practice.
\newline
\newline
$1)$  Suppose all Sylow subgroups of $G$ are cyclic. Then $G \in \mathscr{B}$, $G$ is meta-cyclic or cyclic, $(|[G,G]|,|G/[G,G]|)=1$.
\newline
\newline
$2)$  Suppose $E_8 \lhd G$. Then $G=(E_8 \times R)\langle a,b \rangle$ or $G = (E_8 \times T)\langle c \rangle$ with $R \lhd G$, $R\langle a,b \rangle \in \mathscr{B}_{\pi}$, $R\langle a,b \rangle/R \cong C_7 \rtimes C_3$, $\langle a \rangle \in Syl_7(G)$, $\langle b \rangle \in Syl_3(G)$, $a^7 \in R$, $b^3 \in R$, $E_8 \rtimes \langle a,b \rangle \in \mathscr{B}$, $(2,|R|)=1=(2,|T|)$, $T \lhd G$, $T\langle c \rangle \in \mathscr{B}_{\pi}$, $T\langle c \rangle/T \cong C_7$, $\langle c \rangle \in Syl_7(G)$, $c^7 \in T$, $E_8 \rtimes \langle c \rangle \in \mathscr{B}$. 
\newline
\newline
$3)$ Suppose $E_4 \lhd G$. Then $G = (E_4 \times  U)\langle f \rangle$ with $U \lhd G$, $U\langle f \rangle/U \cong C_3$, $\langle f \rangle \in Syl_3(G)$, $f^3 \in U$, $U\langle f \rangle \in \mathscr{B}_{\pi}$, $(2,|U|)=1$, $E_4\langle f \rangle \in \mathscr{B}$. 
\newline
\newline
$4)$ Suppose $E_{32} < G$. Then $E_{32} \lhd G$ and $G=(E_{32} \times R)\langle a,b \rangle$ with $R\langle a,b \rangle$ as in $2)$, but with the integer $7$ replaced everywhere by $31$ and the integer $3$ replaced everywhere by $5$, $E_{32} \rtimes \langle a,b \rangle \in \mathscr{B}$.  
\newline
\newline
$5)$ Suppose  $Q < G$. Then either $Q \lhd G$ or else $QF(G)/F(G) \lhd G/F(G)$. 
\newline
In the first case one has $G/C_G(Q) \cong A_4$ $\cong(C_2 \times C_2) \rtimes C_3$ and $G/O_{2'}(G) \cong Q \rtimes C_3$ $\cong SL(2,3)$,  with 
$G=(Q \times R)\langle a \rangle$, $\langle a \rangle \in Syl_3(G)$, $a^3 \in R$, $R \lhd G$, $R\langle a \rangle \in \mathscr{B}_{\pi}$, $Q\langle a \rangle \in \mathscr{B}$, $(2,|R|)=1$;
\newline
In the second case one has for $Q \ntrianglelefteq G$  that  $(G/F(G))/(O_{2'}(G/F(G))) \cong$ $Q \rtimes C_3 \cong$ $SL(2,3)$ and that the Sylow $p$-subgroups of $G$ with $p \neq 2$ and $p$ dividing $|G/F(G)|$ are cyclic. Moreover, one has here that $G/F(G) \in \mathscr{B}$ and that either  
$G=((((C_5 \times C_5)^{\delta_1} \times (C_{11} \times C_{11})^{\delta_2}) \rtimes Q) \rtimes T)\langle b \rangle$, $T \lhd G$, $[Q,\langle b \rangle] \neq \{1\}$, $Q \lhd Q\langle b \rangle$,  
$T\langle b \rangle \in \mathscr{B}_{\pi}$, $T\langle b \rangle$ being a Hall-subgroup of $G$,  $\langle b \rangle \in Syl_3(G)$, $b^3 \in T$, $\delta_j \in \{0,1\}$ for $j \in \{1,2\}$ except for $\delta_1=\delta_2=0$, or that $G=(((C_{11} \times C_{11}) \rtimes Q) \times T)\langle b,d \rangle$, $T \lhd G$, $T\langle b,d \rangle$ being a Hall-subgroup of $G$, $\langle b \rangle \in Syl_3(G)$, $b^3 \in T$, $\langle d \rangle \in Syl_5(G)$, $d^5 \in T$, $[Q,\langle b \rangle] \neq \{1\}$, $Q \lhd Q\langle b \rangle$, $[Q,\langle d \rangle]=\{1\}$, $[\langle b \rangle,\langle d \rangle]=\{1\}$, $[C_{11} \times C_{11}, \langle b \rangle] \neq$ $\{1\} \neq [C_{11} \times C_{11},   \langle d \rangle]$. 
\newline
\newline
$6)$ Suppose $C_2 \lhd G$. Then either $G = V\langle c \rangle$, $V \unlhd G$, $(2,|V|)=1$, $\langle c \rangle \in Syl_2(G)$, $G/O_{2'}(G) \cong \langle c \rangle$ or else $Q \in Syl_2(G)$ satisfying $Q \lhd G$. 
\newline
\newline
$7)$ Suppose $E_{p^3} < G$, $p$ an odd prime. Then $E_{p^3} \lhd G$ and $p \nmid |G/E_{p^3}|$. Moreover, as $F(G)=E_{p^3} \times V$ with $V \lhd G$ and $p \nmid |V|$, either it holds that $G/F(G) \in \mathscr{B}$ with all its Sylow subgroups cyclic, or else $G/F(G) \in \mathscr{B}$, where $G/F(G)$ contains a normal Sylow $2$-subgroup isomorphic to $Q$ and where all its other Sylow subgroups are cyclic. Notice $G/E_{p^3} \in \mathscr{B}_{\pi}$ and $G/V \in \mathscr{B}_{\pi}$. 
\newline
\newline
$8)$ Suppose $E_{p^2} < G$, $p$ an odd prime with  $p \nmid |G:E_{p^2}|$. Then $E_{p^2} \lhd G$. Moreover, one has $F(G)=E_{p^2} \rtimes V$ with  $V \lhd G$ and $p \nmid |V|$. Either it holds that $G/F(G) \in \mathscr{B}$ with all its Sylow subgroups cyclic of odd order, or else $p$ divides $55$, $QF(G)/F(G)$ being a normal subgroup of $G/F(G)$ and also satisfying   $G/F(G) \in \mathscr{B}$.   
\newline
\newline
$\textbf{Proof}$ (Sketch) 
\newline 
$\emph{Re 1)}$ See the Corollary to Theorem $3$, $[4]$.
\newline 
$\emph{Re 2) ... Re 5)}$ Look at the relevant results obtained in this article, combined with Gross' work in Theorem $1$, Theorem $2$, Theorem $3$, $[4]$; as to $5)$, see also the Remarks following our Corollary $22$.
\newline
$\emph{Re 6)}$ See the contents of the next Theorem $21$ and Corollary $22$. 
\newline 
$\emph{Re 7)}$ Use the facts that $E_{p^3}$ is isomorphic to the only existing non-cyclic normal subgroup of odd order in $G/V$, and that either no other non-cyclic normal subgroup of prime power order of $G/V$ does exist, or else that also $Q$ is isomorphic to a normal subgroup of $G/F(G)$; look at $5)$ above and also to Theorem $4.1$, $[9]$. 
\newline 
$\emph{Re 8)}$ Similar arguments as in $7)$ do also apply here. $\square$    
\newline
\newline
$\textbf{Theorem 21}$ Let $G \in \mathscr{A}_{\pi}$ and assume its Sylow $2$-subgroups are cyclic. Suppose at least one other Sylow $p$-subgroup of $G$ is not cyclic, call it $A$. Then $A \lhd G$, $A$ is elementary abelian and $F(G)=A \rtimes V$, where $V \unlhd G$ and satisfying $p \nmid |V|$ and $p \nmid |G/A|$. Moreover, $A \cong E_{p^2}$ or $A \cong E_{p^3}$ does hold. In addition, the following structural properties are apparent.
\newline
\newline 
$1)$ Suppose $A \cong E_{p^2}$. Then $G/V \in \mathscr{B}$, whence there exists a $\mathscr{B}$-subgroup $H/V$ of $G/V$ with $H \leq G$, satisfying $|G : H|=p^2$, while $N_{H/V}(C_pV/V)$ is a normal subgroup of $H/V$ for any  $C_p < A$ of order $p$, and  $|(H/V)/N_{H/V}(C_pV/V)|=p+1$. Furthermore, $N_{H/V}(C_pV/V)$ is a meta-cyclic $\mathscr{B}$-group.
\newline
\newline
$2)$ Suppose $A \cong E_{p^3}$. Then $G/V \in \mathscr{B}$. It holds that inside $G$ there are subgroups  $H \geq F(G)$ and  $K \geq F(G)$ such that the quotient group Fitting group $(G/V)/F(G/V)$ is a direct product of the groups $(H/V)/F(G/V)$ and $(K/V)/F(G/V)$  (say), where $|(H/V)/F(G/V)|$ is odd and where 
$(H/V)/F(G/V)$ has all its Sylow subgroups cyclic, $H/V \in \mathscr{B}$, and with $(K/V)/F(G/V)$ being a cyclic $2$-group normalizing each subgroup of $F(G/V)$.
\newline
\newline
$\textbf{Proof}$ All the non-cyclic Sylow subgroups of $G$ are contained in $F(G)$; note that $G$ is solvable, as $G$ is $2$-nilpotent due to Burnside's Theorem III.$2.8$ Satz, $[5]$, because it is given that the Sylow $2$-subgroups of $G$ are cyclic. Therefore, as $G \in \mathscr{A}_{\pi}$, $F(G)$ is a direct product of $A \lhd G$  and  $V \lhd G$ satisfying $(|A|,|V|)=1$, where $A$ is some non-cyclic elementary abelian $p$-subgroup of $G$ for some specific odd prime $p$, as assumed in the statement of the Theorem. We know already that now $A \cong E_{p^2}$ or $A \cong E_{p^3}$ has to be the case. 
\newline
$\emph{Re 1)}$ Assume $A \cong E_{p^2}$. Then $G/V$ does contain precisely one non-cyclic Sylow subgroup (namely $AV/V$), whence $AV/V \lhd G/V$ does hold. Hence from Gross' Theorem $3$, $[4]$ one deduces $G/V \in \mathscr{B}$. The rest of the statements in the Theorem under $1)$ is now nothing else but a recasting of the contents of Theorem $6.1$, $[10]$.
\newline
$\emph{Re 2)}$ Assume $A \cong E_{p^3}$. Here too one gets $G/V \in \mathscr{B}$ from Gross' Theorem $3$, $[4]$, as $AV/V$ happens to be the only existing non-cyclic Sylow subgroup of $G/V$. Now this time the rest of the statements of the Theorem under $2)$ is a recasting of the contents of Theorem $8.1$ of $[10]$. 
\newline
The proof of the Theorem is complete. $\square$
\newline
\newline
$\textbf{Corollary 22}$ Let $G \in \mathscr{A}_{\pi}$ be solvable. Assume that $Q$ is not isomorphic to a normal subgroup of $G$ and that $G$ contains a normal subgroup of order $2$. Then either $G$ is a $\mathscr{B}$- group with all its Sylow subgroups cyclic, or else $G$ is a group dealt with in Theorem $21$. In both cases, $G$ is $2$-nilpotent whose Sylow $2$-subgroups are cyclic; compare with $6)$ in Theorem $20$. 
\newline
\newline
$\textbf{Proof}$ The given group $G \in \mathscr{A}_{\pi}$ admits by assumption precisely one element $a$ (say) of order $2$. Hence either any Sylow $2$-subgroup of $G$ is cyclic or generalized quaternion by III.$8.2$ Satz, $[5]$. If generalized quaternion, it has to be of order $8$, because of $G \in \mathscr{A}_{\pi}$; call it $Q$. In that case, it holds that some Sylow $2$-subgroup $S$ of $G/\langle a \rangle$ is isomorphic to $E_4$. In fact, that subgroup $S$ is the only existing subgroup of order $4$ in $G/\langle a \rangle$; see $(4.4)$ Theorem, $[7]$. Hence $Q$ would be a normal subgroup of $G$, contrary to the assumptions made in the Corollary. Thus the statement in the Corollary does hold by what has been put in the assumptions of Theorem $21$.  $\square$
\newline
\newline
$\textbf{Remark}$ In Theorem $20$, $5)$ it was investigated what happens when $Q \leq G$ for a solvable $\mathscr{A}_{\pi}$-group $G$. Here we inspect carefully where $Q \ntrianglelefteq G$ leads to for such groups $G$. Notice $Q \in Syl_2(G)$. It holds that $Q$ intersected with the Fitting subgroup $F(G)$ of $G$ provides $Q \cap F(G)=\{1\}$.
\newline
\newline
[Indeed, if $Q \cap F(G) \neq \{1\}$, then     $Z(Q) \leq F(G)$ and because of $Q \ntrianglelefteq G$, $F(G)$ does contain a normal Sylow $2$-subgroup $S$ (say) of order at most $4$ satisfying  $S \geq Z(Q)$; remember that $F(G)$ is nilpotent, that each Sylow $2$-subgroup of $G$ is isomorphic to $Q$, that $Q$ contains as subgroups of order $4$ only cyclic ones, and that $Z(Q)$ is the only existing subgroup of order $2$ inside $Q$. So, if $|S|=4$ would hold, then, by $G \in \mathscr{A}_{\pi}$, any cyclic subgroup of $G$ being of order $4$ would be contained in $F(G)$, whence in particular it happens that $Q \leq F(G)$, as $Q$ is generated by its cyclic subgroups of order $4$. By assumption though it is therefore not possible that $|S| = 4$. 
\newline
Thus $Z(Q) \in Syl_2(F(G))$ satisfying $Z(Q) \unlhd F(G)$. Then, however, $Z(Q)$ is characteristic in $F(G)$, yielding $Z(Q) \unlhd G$. Thus $G/Z(Q) \in \mathscr{A}_{\pi}$, where $Q/Z(Q)$ is isomorphic to $E_4$ and where now  $Q/Z(Q) \in Syl_2(G/Z(Q))$. Hence $Q/Z(Q) \unlhd G/Z(Q)$ by $(4.4)$ Theorem, $[7]$. Hence $Q$ would be a normal subgroup of $G$ yielding $Q \unlhd F(G)$, which is not allowed.]
\newline
\newline
Thus  $QF(G)/F(G)$ is a Sylow $2$-subgroup of $G/F(G)$ being isomorphic to $Q$. Remember, because of $G \in \mathscr{A}_{\pi}$, that all non-cyclic elementary abelian subgroups of $G$ are contained in $F(G)$. Any other Sylow subgroup of $F(G)$ is thus cyclic as there are no Sylow subgroups of $F(G)$ existent being non-cyclic and simultaneously not elementary abelian. Hence each Sylow subgroup of $G/F(G)$ is cyclic in case it is not isomorphic to $Q$. Recall that $G/F(G) \in \mathscr{A}_{\pi}$. Then, just as it has been argued above, either $QF(G)/F(G)$ is normal in $F(G/F(G))$ or otherwise
$$QF(G)/F(G) \cap F(G/F(G))=\{1\} \hspace{3cm} (\star)$$

satisfying $QF(G)/F(G) \cong Q$. 
\newline
\newline 
Suppose the last equality $(\star)$ is in vogue. Let us introduce the notations $\overline{G}=G/F(G)$, $\overline{Q}= QF(G)/F(G)$ and $\overline{F}=F(G/F(G))$.  It holds that all the Sylow subgroups of $\overline{F}$ are cyclic. Since $\overline{F}$ is therefore abelian, it holds that not only $2 \nmid |\overline{F}|$ is fulfilled, but also  that $\overline{F}$ is the direct product of cyclic subgroups of odd order being pairwise relatively prime. Since $\overline{G}$ is solvable, it is a fact that $C_{\overline{G}}(\overline{F})=\overline{F}$; 
see III.$4.2$ Satz.$[5]$. Hence by the structure of $\overline{F}$ obtained above, one gets that $\overline{G}/C_{\overline{G}}(\overline{F})$ is abelian; notice, that group is a subgroup of $Aut(\overline{F})$ and that here $Aut(\overline{F})$ is abelian. We have thus obtained a contradiction to $(*)$.
\newline
\newline
In summary, if $Q \leq G$, $G \in  \mathscr{A}_{\pi}$ being solvable, then either $Q \unlhd G$ or otherwise $\overline{Q}$ is normal in $\overline{G}$ with $\overline{Q} \cong Q$; in the last situation $G$ is a      $\mathscr{B}$-group as it satisfies the assumptions of Theorem $11$, $[9]$, whose structure has been fully described in Theorem $4.1$, $[10]$. Knowing that the group $\overline{G}$ with $\overline{Q} \lhd \overline{G}$ contains only as other Sylow subgroups cyclic ones, one therefore gets, that there exists a $3$-nilpotent meta-cyclic       $\mathscr{B}$-subgroup $L\langle a \rangle$ of $\overline{G}$ such that $\overline{G}=(\overline{Q} \times L)\langle a \rangle$ in which $L$  is a meta-cyclic      $\mathscr{B}$-group satisfying $\langle a \rangle \in Syl_3(\overline{G})$, $[\langle a \rangle,\overline{Q}] \neq \{1\}$, $a^3 \in L$, $L \unlhd L\langle a \rangle$, $\overline{Q}\langle a \rangle/\langle a^3 \rangle \cong SL(2,3)$. As such these facts do sharpen the structure of the groups $T\langle b \rangle$ and $T\langle b,d \rangle$ mentioned in Theorem $20$, $5)$.
\section{References}
\bibliography{sample}
$[1]$ E. Artin, The orders of the classical simple groups, Comm. in Pure and Appl. Math., $8 \text{ } (1955), 455-472$.
\newline 
$[2]$ M. Costantini and E. Jabara, On finite groups in which cyclic subgroups of the same order are conjugate, Comm. in Alg., $37 \text{ } (2009), 3966-3990$.
\newline 
$[3]$ A. Bensaïd and R.W. van der Waall, Non-solvable finite groups whose subgroups of equal order are conjugate, Indag. Math., New Series, $1 \text{ }(1990), 397-408$.
\newline 
$[4]$ F. Gross, Finite groups in which any two primary subgroups of the same order are conjugate, Michigan Math. J. $19 \text{ } (1972), 333-339$.
\newline 
$[5]$ B. Huppert, Endliche Gruppen I, Grundlagen der Math. Wissenschaften, $134$, Springer Verlag, Berlin, $1967$.
\newline 
$[5a]$ B. Huppert and N. Blackburn, Finite Groups II, Grundlagen der Math. Wissenschaften, $242$, Springer Verlag, Berlin, $1982$.
\newline
$[5b]$ B. Huppert and N. Blackburn, Finite Groups III, Grundlagen der Math. Wissenschaften, $243$, Springer Verlag, Berlin, $1983$.
\newline
$[6]$ R.C. Lindenbergh, R.W. van der Waall, Ergebnisse über Dedekind-Zeta-Funktionen, monomiale Charaktere und Konjugationsklassen endlicher Gruppen, unter Benutzung von GAP, Bayreuther Math. Schriften, Heft $56$ $(1999), 79-148$.
\newline
$[7]$ S. Sezer, On finite solvable groups all of whose cyclic $p$-subgroups of equal order are conjugate, J. of Alg., $415 \text{ } (2014), 214-233$.
\newline
$[8]$ S. Sezer and R.W. van der Waall, Finite groups all of whose abelian subgroups of equal order are conjugate, Turkish J. of Math., $30 \text{ } (2006), 139-175$.
\newline
$[9]$ R.W. van der Waall, Finite groups whose subgroups of equal order are conjugate, Indag. Math., New Series, $4 \text{ } (1993), 239-254$.
\newline
$[10]$ R.W. van der Waall, The classification of the finite groups whose subgroups of equal order are conjugate, Indag. Math., New Series, $23 \text{ } (2012), 448-478$ and Indag. Math., New Series, $24 \text{ } (2012), 489-490$.
\newline
$[11]$ R.W. van der Waall, The classification of the finite groups whose supersolvable (nilpotent) subgroups of equal order are conjugate, Indag. Math., New Series, $26 \text{ } (2015), 380-383$. 
\newline
$[12]$ R.W. van der Waall and S. Sezer, On finite non-solvable groups whose cyclic $p$-subgroups of equal order are conjugate, submitted.

\end{flushleft}
\end{document}